\documentclass[12pt,a4paper]{amsart}
\usepackage{epsfig}
\usepackage{amssymb}
\usepackage{amsmath}
\usepackage{latexsym}
\usepackage{color}
\usepackage{graphicx}
\usepackage{float}
\usepackage{amsthm}
\usepackage{fullpage}
\newtheorem{definition}{Definition}[section]
  
   \newtheorem{lemma}[definition]{Lemma}
  \newtheorem{corollary}[definition]{Corollary}

\date{\today}
\begin{document}
\title{A remark on Whitehead's cut-vertex Lemma} 
\author{Michael Heusener and Richard Weidmann}

\maketitle

\begin{abstract} We observe that Whitehead's cut-vertex lemma is an immediate consequence of Stallings folds.
\end{abstract}

\section*{Introduction}  Whitehead's cut-vertex lemma \cite{WH} implies that the Whitehead graph of a primitive element of the free group is either not connected or has a {cut vertex}. This result was generalized by Stong~\cite{Stong} and Stallings \cite{St}  to separable subsets of free groups. In this note we observe that this is an essentially trivial consequence of Stallings folds. To do so we relate the connectivity properties of the Whitehead graph of a set of conjugacy classes of elements to the readability of cyclically reduced representatives of these classes in a class of labeled graphs that we call almost-roses.

\section{Free groups and graphs}  

 In the following we denote the free group in $X=\{x_1,\ldots ,x_n\}$ by $F_n$. Throughout this note $n\ge 2$ and $X$ are fixed. As usual we represent elements of $F_n$ by words in $\tilde X:=X\cup X^{-1}$. For any $g\in F_n$ we denote the conjugacy class of $g$ by $[g]$ and for $S\subset F_n$ we set $[S]:=\{[s]\mid s\in S\}$. We say that $g\in F_n$ is primitive if $g$ in contained in a basis of $F_n$. We further call a set $S\subset F_n$ separable if there exists a non-trivial decomposition $F_n=F^1*F^2$ such that any $w\in S$ is conjugate into $F^1$ or $F^2$, see~\cite{St}. 
 Clearly $\{g\}$ is separable for any primitive element $g\in F_n$.
 
 \smallskip  An $\tilde X$-labeled graph $\Gamma$ is a directed graph  (with inverse edges)  with a labeling map $\ell:E\Gamma\to \tilde X$ such that $\ell(e^{-1})=\ell(e)^{-1}$ for any edge $e\in E\Gamma$.  Multiple edges and loops are permitted. The label of a path 
 $\gamma = e_1\cdots e_k$ in $\Gamma$ is
the word $\ell(\gamma)=\ell(e_1)\cdots\ell(e_k)$.
 
\smallskip In the following $R_n$ is the $\tilde X$-labeled graph whose underlying graph has vertex set $V=\{v_0\}$ and edge set $\{e_1^{\pm 1},\ldots ,e_n^{\pm 1}\}$ such that $\ell(e_i^{\varepsilon})=x_i^{\varepsilon}$ for all $\varepsilon\in\{-1,1\}$ and $1\le i\le n$. As usual we identify $F_n$ with $\pi_1(R_n,v_0)$ in the obvious way. For any $\tilde X$-labeled graph $\Gamma$ there is a unique label preserving morphism $f:\Gamma\to R_n$. 
 
\smallskip We say that a word $w$ is readable in $\Gamma$ if there exists a closed path  
%
%
in $\Gamma$ with label $w$.
We will say that an element $g\in F_n$ is readable in $\Gamma$ if the corresponding reduced word 
is readable in $\Gamma$. We will further say that a conjugacy class $[g]$  is readable in $\Gamma$ if some (and therefore all) cyclically reduced word representating  $[g]$  is readable in $\Gamma$.

We will be mostly interested in connected  $\tilde X$-labeled graphs $\Gamma$ such that associated map $f_*:\pi_1(\Gamma)\to \pi_1(R_n)=F_n$ is surjective. In this case the morphism $f$ can be written as a product of Stallings folds, see \cite{St0}. 

\smallskip Recall that a graph $\Gamma$ is called a core graph  if it contains no proper subgraph such that the inclusion is a homotopy equivalence; for a finite graph this just means that $\Gamma$ has no vertex of valence $1$. We also call a pair $(\Gamma,v_0)$ a core pair if $\Gamma$ contains no proper subgraph containing the vertex $v_0$ such that the inclusion is a homotopy equivalence, for a finite graph this means that it contains no vertex of valence $1$ distinct from $v_0$.

\section{Tame sets of elements in the free group}

We consider words that can be read in a simple class of labeled graphs which we call almost-roses.

\begin{definition} We call a $\tilde X$-labeled core graph $\Theta$ of Betti number $n$ an {\em almost-rose} if it folds onto $R_n$ with a single fold. 
\end{definition}

An almost-rose $\Theta$  must have $n+1$ edges and $2$ vertices $u$ and $v$. The fold must identify $u$ and $v$ and preserve the Betti number, it follows that it must fold a loop edge to a non-loop edge. Thus after permuting and/or inverting the $x_i$ and possibly exchanging $u$ and $v$ there exist $1\le k\le l\le n$ with $k<n$ such that the following hold:
\begin{enumerate}
\item It has a loop edge based at $u$ and an edge from $u$ to $v$ labeled by $x_1$.
\item There are loop edges based at $u$ with labels $x_2,\ldots,x_k$.
\item There are edges from $u$ to $v$ with labels $x_{k+1},\ldots x_l$.
\item There are loop edges based at $v$ with lables $x_{l+1},\ldots ,x_n$.
\end{enumerate}
We denote this almost rose by $\Theta_{k,l}^{n}$. 

\begin{figure}[htb]
  \centering 
  \def\svgwidth{225pt}
 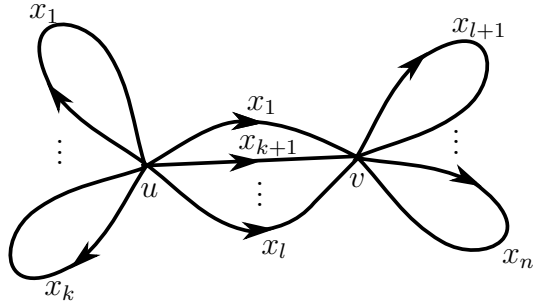
  \caption{The almost-rose $\Theta_{k,l}^{n}$.}\label{almostrose}
\end{figure}

\begin{definition} We call a set of conjugacy classes $[S]$ of elements of  $F_n$ {\em tame} if there exists  an almost-rose $\Theta$ such that any $[s]\in [S]$ is readable in $\Theta$. We further call $g$ tame if $\{[g]\}$ is tame.
\end{definition}


\begin{lemma}\label{primitiveimpliestame} If $g\in F_n$ is  primitive  then $[g]$ is tame.
\end{lemma}

\begin{proof} Let $g\in F_n$ be primitive and $w$ be a cyclically reduced word representing $[g]$. If $w\in \tilde X$ then the claim is trivial. Thus we may assume that $w$ is  a cyclically  reduced word of length at least $2$. Choose a tuple of reduced words $(w_1=w,w_2,\ldots ,w_n)$ representing a basis of $F_n$. Let $\Gamma$ be the wedge of $n$ (subdivided)  circles joined at a vertex $u_0$ such that the path travelling the $i$-th circle starting and ending at $u_0$ is labeled by $w_i$. Note that $\Gamma\neq R_n$ as $w_1$ is of length at least $2$.

The associated morphism $f:\Gamma\to R_n$ is $\pi_1$-surjective and therefore folds onto $R_n$. As $w=w_1$ is cyclically reduced and as all $w_i$ are reduced it follows that any vertex is the  terminal vertex of two edges with different labels. It follows that any graph occuring in the folding sequence is a core graph, thus $f$ factors through an almost-rose $\Theta$. As any word readable in $\Gamma$  is also readable in $\Theta$ { and as $w=w_1$ is readable in $\Gamma$ it follows that $w$} is readable in $\Theta$. Thus $[g]$ is tame.
\end{proof}


The argument easily generalizes to separable sets.

\begin{lemma}\label{separableimpliestame} If $S\subset F_n$ is separable  then $[S]$ is tame.
\end{lemma}

\begin{proof} Choose a free decomposition $F_n=F^1*F^2$ such that any $s\in S$ is conjugate to an element of $F^1$ or $F^2$. Possibly after conjugating $F^1$ and $F^2$ with some $g\in F_n-1$ we may assume that there exists no $x\in\tilde X$ such that any $h\in (F^1\cup F^2)-\{1\}$ can be represented by a reduced word of type $xwx^{-1}$.

For $i=1,2$ choose a labeled core pair $(\Gamma_i,v_i)$ such that any element of $F^i$ can be read by a closed path in $\Gamma_i$ based at $v_i$. The existence of such a pair follows from the theory of covering spaces, alternatively it can be constructed using folds starting with a wedge of circles labeled by a generating set of $F^i$. It follows that $[s]$ is readable in $\Gamma_1$ or $\Gamma_2$  for any $s\in S$.

We now construct a new labeled graph $\Gamma$ by taking the wedge of $\Gamma_1$ and $\Gamma_2$ at $v_1=v_2$. The above assumption on the free decomposition guarantees that any vertex is the terminal vertex of two edges with different labels. Note that the label preserving map $f:\Gamma\to R_n$ is $\pi_1$-surjective as $\langle F^1,F^2\rangle=F_n$. We distinguish two cases.

If $f$ is an isomorphism then $F^1$ and $F^2$ are free factors generated by subsets of the basis $X$. Then we can find an almost-rose $\Theta$ (with $k=l$) such that $[s]$ can be read in $\Theta$ for any $s\in S$; indeed simply take $\Theta$ to be the graph obtained from $\Gamma_1$ and $\Gamma_2$ joined by an edge with arbitrary label. Thus $S$ is tame.

If $f$ is not an isomorphism then we see as in the proof of the previous lemma that $f$ factors through an almost-rose $\Theta$. It follows that $[s]$ is readable in $\Theta$ for any $s\in S$ as it readable in $\Gamma$. It follows that $[S]$ is tame.
\end{proof}

\section{Whitehead graphs}

We define the Whitehead graph $\mathrm{Wh}(\Gamma)$ of a $\tilde X$-labeled graph $\Gamma$ to be the (undirected) graph whose vertex set is $\{x_1,x_1^{-1},\ldots ,x_n,x_n^{-1}\}$ such that distinct vertices $x$ and $y$ are joined by an edge if and only if the word $xy^{-1}$ 
occurs as a label in a reduced path in $\Gamma$.
\begin{figure}[htb]
\centering
 \def\svgwidth{\columnwidth}\footnotesize
 \begin{center}
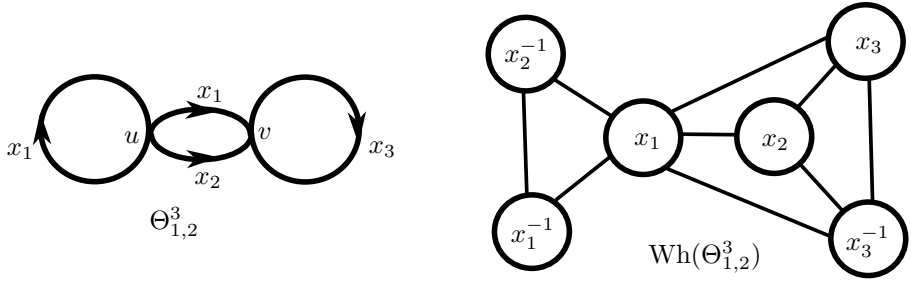 
 \end{center}
  \caption{The almost-rose $\Theta^{3}_{1,2}$ and its Whitehead graph.}\label{exampleWh}
\end{figure}

 {Recall that a vertex $v$ of a connected graph $\Gamma$ is called a cut vertex if its removal disconnects the component of $\Gamma$ containing $v$.}

\begin{lemma}\label{connectivityofwg} Let $\Theta$ be an almost-rose. Then $\mathrm{Wh}(\Theta)$ has a {cut vertex}. 
\end{lemma}

\begin{proof} Without loss of generality we may assume that $\Theta=\Theta^n_{k,l}$ for some $k$ and $l$.
With the notation from the previous section it follows that  $\mathrm{Wh}(\Theta^n_{k,l})$ is the wedge at $x_1$ of the complete graphs with vertex sets 
\[
V_1:=\{x_1^{\pm 1},\ldots ,x_k^{\pm 1},x_{k+1}^{-1},\ldots ,x_l^{-1}\}
\ \text{ and }\ V_2:=\{x_1,x_{k+1},\ldots x_l,x_{l+1}^{\pm 1},\ldots ,x_n^{\pm 1}\}
\]
(see Figure~\ref{exampleWh}). Thus $x_1$ is a {cut vertex} of $\mathrm{Wh}(\Theta)$.
\end{proof}

Clearly the existence of a label preserving morphism $f:\Gamma\to\Gamma'$ implies that $\mathrm{Wh}(\Gamma)\subset \mathrm{Wh}(\Gamma')$. If $\Gamma'$ is an almost rose then the converse holds as well:

\begin{lemma}\label{inclusionofwhiteheadgraphs} Let $\Gamma$ be an $\tilde X$-labeled graph and $\Theta$ an almost rose. Then the following are equivalent:
\begin{enumerate}
\item There exists a label preserving map $f:\Gamma\to \Theta$.
\item $\mathrm{Wh}(\Gamma)\subset \mathrm{Wh}(\Theta)$.
\end{enumerate}
\end{lemma}
\begin{proof} Note that (1) implies (2) by the remark preceding this lemma.

Suppose that $\mathrm{Wh}(\Gamma)\subset \mathrm{Wh}(\Theta)$, we may assume that $\Theta=\Theta^n_{k,l}$ for some $k,l$. We define a map $f_E:E\Gamma\to E\Theta$ as follows:
\begin{enumerate}
\item If $\ell(e)=x_j\neq x_1$ then $f_E(e)$ is defined to be the unique edge of $\Theta$ labeled by $x_j$. 
\item If $\ell(e)=x_1$ and there exists an edge $e'\in E\Gamma$ with the same terminal vertex as $e$ such that $\ell(e')\in V_2$ (Here $V_2$ is defined as in the proof of Lemma~\ref{connectivityofwg}) we define $f_E(e)$ to be the edge of $\Theta$ between $u$ and $v$ that has label $x_1$. 
\item Otherwise we define $f_E(e)$ to be the loop edge of $\Theta$ with label $x_1$. 
\end{enumerate}
It is easily verified that $\mathrm{Wh}(\Gamma)\subset \mathrm{Wh}(\Theta)$ guarantees that $f_E$ extends to a graph morphism $f:\Gamma\to\Theta$.
\end{proof}

For a set $S\subset F_n$ the Whitehead graph $\mathrm{Wh}(S)$ is defined as the graph with vertex set $\{x_1,x_1^{-1},\ldots ,x_n,x_n^{-1}\}$ such that vertices $x$ and $y$ are joined by an edge if the word $xy^{-1}$ occurs as a subword of a {power of a} cyclically reduced conjugate of some element of $S$. For any $g\in F_n$ we further define $\mathrm{Wh}(g):=\mathrm{Wh}(\{g\})$. 

\smallskip For any $S\subset F_n$ we define $\Gamma_S$ to be the (not connected) $\tilde X$-labeled graph that is the disjoint union of circuits whose labels are cyclically reduced words representing the conjugacy class of $s$ for each $s\in S$. 
It is obvious that $\mathrm{Wh}(\Gamma_S)=\mathrm{Wh}(S)$.

\begin{corollary}\label{connectivityofwg2} Let $S\subset F_n$. Then the following are equivalent:
\begin{enumerate}
\item $[S]$ is tame.
\item  $\mathrm{Wh}(S)$ is not connected or $\mathrm{Wh}(S)$ has a {cut vertex}.
\end{enumerate}
\end{corollary}

\begin{proof} If $[S]$ is tame then there exists an almost-rose $\Theta$ such that for any $[s]\in [S]$ any cyclically reduced representative of $[s]$ is readable in $\Theta$, i.e. there is a label-preserving map from $\Gamma_S$ to $\Theta$. Thus $$\mathrm{Wh}(S)=\mathrm{Wh}(\Gamma_S)\subset \mathrm{Wh}(\Theta).$$ As $\mathrm{Wh}(\Theta)$ has a {cut vertex} by Lemma~\ref{connectivityofwg} it follows that $\mathrm{Wh}(S)$ is either not connected or has a {cut vertex} as both graphs have the same vertex sets.

If $\mathrm{Wh}(\Gamma_S)$ is not conntected or has a cut vertex then $\mathrm{Wh}(S)=\mathrm{Wh}(\Gamma_S)\subset \mathrm{Wh}(\Theta)$ for some almost rose $\Theta$ which implies that there exists a label preserving morphism $f:\Gamma_S\to\Theta$ which implies that $[S]$ is tame.\end{proof}

The following corollaries are immediate consequences of Corollary~\ref{connectivityofwg2}, 
Lemma~\ref{primitiveimpliestame} and Lemma~\ref{separableimpliestame}, respectively.

\begin{corollary}[Whitehead]\label{corwhitehead} Let $w\in F_n$ be primitive. Then $\mathrm{Wh}(w)$ is not connected or has a {cut vertex}.
\end{corollary}

\begin{corollary}[Stallings, Stong]\label{corstallings} Let $S$ be a separable set. Then $\mathrm{Wh}(S)$ is not connected or has a {cut vertex}.
\end{corollary}

\noindent{\em Remark: } Note that in order to prove Corollary~\ref{corwhitehead} and Corollary~\ref{corstallings} we only need the trivial implication of Lemma~\ref{inclusionofwhiteheadgraphs} and the first implication of Lemma~\ref{connectivityofwg}.

\section{Other variants of the Whitehead graph}
As has been pointed out to the authors by Warren Dicks, the graph we define is not the graph used in the original paper of Whitehead, see \cite{D} for a detailed discussion. In particular Whitehead deals with both conjugacy classes of elements and elements. However the connectivity properties of the graphs considered by Whitehead can also be studied using the ideas presented in this note.

Recently, Henry Wilton also gave a proof of Whitehead's lemma using foldings. See Lemma 2.10 in \cite{Wilton}.

\end{document}